\documentclass[12pt,oneside]{article}
\usepackage{amsmath,amsthm,amsfonts,amssymb,MnSymbol,mathrsfs,latexsym}
\usepackage{mathtools}
\usepackage[english]{babel}

\newcommand{\const}{\rm const}

  \newcommand{\Dom}{\rm  Dom}

\DeclareMathOperator*{\esssup}{ess\,sup}

\theoremstyle{plain}
\newtheorem{theorem}{Theorem}[section]

\newtheorem{lemma}[theorem]{Lemma}
\newtheorem{proposition}[theorem]{Proposition}

\newtheorem{remark}{Remark}[section]
\newtheorem{example}{Example}[section]

\renewenvironment{proof}{{\bf{Proof.}}}{\hfill $\Box$ \\}

\title{\large \textbf{Factorizable convergence of random variables\\ in Grand Lebesgue Spaces}}

\footnotesize\date{}

\author{\normalsize Maria Rosaria Formica ${}^{1}$,   \normalsize Eugeny Ostrovsky
${}^2$ and \normalsize Leonid Sirota ${}^3$}

\begin{document}

  \maketitle

\begin{center}
{\footnotesize ${}^{1}$ Universit\`{a} degli Studi di Napoli \lq\lq Parthenope\rq\rq, via Generale Parisi 13,\\
Palazzo Pacanowsky, 80132,
Napoli, Italy.} \\

\vspace{1mm}

{\footnotesize e-mail: mara.formica@uniparthenope.it} \\

\vspace{2mm}

{\footnotesize ${}^{2,\, 3}$  Bar-Ilan University, Department of Mathematics and Statistics, \\
52900, Ramat Gan, Israel.} \\

\vspace{1mm}

{\footnotesize e-mail: eugostrovsky@list.ru}\\

\vspace{1mm}

{\footnotesize e-mail: sirota3@bezeqint.net} \\

\end{center}

\vspace{3mm}

\begin{abstract}
 \hspace{3mm} We obtain results concerning the so-called {\it factorization} for the convergence of random variables {\it almost everywhere} (almost surely or with probability one),
  belonging to the classical Lebesgue-Riesz spaces and we extend these results to the Grand Lebesgue Spaces.\par
  \ We also give exact estimates for the parameters involved and provide several examples. \par
 \ We also show that the obtained estimates are, in the general case, essentially non-improvable, of course up to a multiplicative constant.\par
\end{abstract}



\vspace{3mm}


\vspace{3mm}


 \vspace{5mm}

\section{Introduction.}

\vspace{1mm}

P.E. Kloeden and A. Neuenkirch in \cite{Kloeden Neuekirch} proved a result which establishes a connection between the convergence rates in the p-th mean ($p\geq 1$) and the pathwise convergence rates for sequences of random variables belonging to the classical Lebesgue-Riesz spaces. The importance of this result is also due to applications, given by the authors, for  It\^{o} stochastic differential equations.

In this paper we extend the above result to random variables belonging to the generalized Grand Lebesgue Spaces. These spaces have a crucial role in several applications, as in Probability theory, Functional Analysis, interpolation theory, PDEs.

\vspace{1mm}

 \hspace{3mm} Let $ (\Omega = \{\omega\}, \,  {\bf B}, \, {\bf P} ) $ be a probability space. Let $ \{\xi_n \}, \ n =1, 2,3,\ldots$, be a sequence of numerical valued random variables (r.v.). We say that $\{\xi_n \}$ converges to zero as $n \to \infty $ {\it almost everywhere (a.e.)}, or equally {\it almost surely (a.s)}, or
{\it with probability one (w.p.o.)} if
\begin{equation} \label{gonerg ae}
{\bf P} \left(\lim_{n \to \infty} \xi_n = 0 \ \right) = 1.
\end{equation}

The following theorem can be found, e.g., in \cite{Kantorovich Akilov}, chapter 1, section 6, Theorem 5, part 2, known as \emph{Theorem of regulator} of the convergence a.e. or, in the other words, {\it factorizable convergence.}

\begin{theorem}\label{theorem regulator} {\rm \textbf{(Regulator)}}
Let $ \{\xi_n \}, \ n =1, 2,3,\ldots$, be a sequence of numerical valued random variables converging to zero {\it almost everywhere (a.e.)}. Then there exists a non-negative finite a.e. random variable $ \ \zeta \ $ and a sequence of {\it non-random} real valued non-negative variables $ \epsilon_n $ tending to  zero ($ \epsilon_n \to 0 $) such that
\begin{equation} \label{regulator}
|\xi_n| \le \epsilon_n \ \zeta.
\end{equation}

\end{theorem}

\ Of course, the converse implication is also true.

 \vspace{3mm}

 The pair $ \  \{ \epsilon_n, \ \zeta\}, \ $  non unique, may be named by definition as the {\it \textbf{regulator}} for the sequence $ \ \{\xi_n \}. \ $

\vspace{3mm}

 \begin{remark}
 {\rm
 Note that the famous book \cite{Kantorovich Akilov} does not contain the proof of this assertion as well as the references of the works where this proposition is ground. For the filling of this gap see \cite{Ostrovsky5,Ostrovsky7}. For instance, \eqref{gonerg ae} holds true if and only if
\begin{equation} \label{criterion}
 \lim_{n \to \infty} {\bf E} \sup_{m \ge n} [ \ |\xi_m|/(1 + |\xi_m|) \ ] =  0.
\end{equation}

\vspace{1mm}

\hspace{3mm} Really, if the relation \eqref{gonerg ae} is satisfied, then \eqref{criterion} is grounded on the basis of  the dominated convergence theorem.\par
 \ Conversely, let the condition \eqref{criterion}  be satisfied. There exists a deterministic positive sequence $ \,\{\delta_n \} \ $ for which $\displaystyle \lim_{n \to \infty} \delta_n = 0 \ $ and

$$
\lim_{n \to \infty} {\bf E} \sup_{m \ge n} [ \ |\xi_m/\delta_m|/(1 + |\xi_m|/\delta_m) \ ] =  0.
$$
Therefore $ \displaystyle \lim_{n \to \infty}[ \ |\xi_n|/\delta_n \ ] = 0, \ $ with probability one, and denoting
 $ \ \upsilon := \displaystyle\sup_n [ \ |\xi_n|/\delta_n \ ] \ $ we have
 $$
 {\bf P} (0 \le \upsilon < \infty) = 1
 $$
 and finally
 $$
 |\xi_n| \le \upsilon \cdot \delta_n,
 $$
i.e. the existence of a regulator of the convergence.
}
 \end{remark}

 \vspace{4mm}

 \ An useful estimation for this functional, in the case when the r.v. $ \  \{ \xi_n  \} \ $ are sums of dependent variables, is offered in the recent article \cite{Da Silva}. \par

\vspace{3mm}

  Moreover, inequality \eqref{regulator} remains true if the r.v. $ \{\xi_n\} $ take values in a {\it separable} Banach space $ B = (B,||\cdot||_B), \ $
 where the absolute value $ |\cdot |  $ in the inequality is replaced by the norm $ ||\cdot||_B $ (see \cite{Ostrovsky4,Ostrovsky5}).

\vspace{3mm}

\begin{remark}\label{another criterion}
{\rm
 Another criterion for the convergence a.e. (a.s.) states that a sequence of numerical valued r.v. $ \ \{X_n\}, \ n = 1,2, \ldots$, converges a.e. to zero iff
\begin{equation} \label{iff}
 \forall \epsilon > 0 \ \Rightarrow  \ \lim_{n \to \infty} {\bf P} \left( \bigcup_{m \ge n}  \{ |X_m| \le \epsilon \}  \ \right) = 1.
 \end{equation}

 }
 \end{remark}

\hspace{3mm} But our approach is more convenient for us. \par

 \vspace{3mm}

\begin{remark}
{\rm
Let $ \{X_n\}$ be a sequence of numerical valued r.v. and consider the sequence
$$
 \eta_n :=  \frac{|X_n|}{1 + |X_n|},
$$
which is positive and bounded, $  0 \le \eta_n \le 1 $.

 The almost everywhere convergence to zero of $  \{X_n\}, \ n \to \infty$, is completely equivalent to the following convergence
 $$
 \lim_{n \to \infty} {\bf E} \left\{ \ \sup_{m \ge n} \eta_m \  \right\} = 0.
 $$

\vspace{3mm}

Note that the  metric (symmetry, positivity, triangle inequality) and other properties of the distance function
$$
\rho(X,Y) := {\bf E} \frac{|X-Y|}{1 + |X - Y|}
$$
in more than used in this report (one-dimensional case  $ \ \dim(X) = \dim(Y) = 1 $)
are in detail investigated in the recent articles \cite{Amaro Bessa Vilarinho,Arbieto Bochi}.
}
\end{remark}

\vspace{3mm}

\begin{remark}{\it Continuum case.}
{\rm
Let now $ \ \eta(t), \ t \in [0,1]$, be a separable bi - measurable {\it random process} such
 that $ \ \eta(1) = 0. \ $ \textbf{Question}: under what conditions imposed on the  distribution of the random process $ \ \eta(t), \ t \in [0,1]$,
\begin{equation} \label{what conditions}
{\bf P}(\lim_{t \to 1^-} \eta(t) = 0) = 1  \ ?
\end{equation}
}
\end{remark}

\vspace{3mm}

  Note first of all that there exists a simple (but not convenient in practice)
  {\it necessary and sufficient}  condition for \eqref{what conditions}: indeed, when there exists a constant
 $ \ \delta \in (0,1) \ $ such that almost all trajectories of the random process $ \ \eta(\cdot) \ $ are bounded in the interval $ \ [1 - \delta, 1] $, i.e.
\begin{equation} \label{first cond}
\exists \,\delta \in (0,1) \ : \ {\bf P}(\sup_{t \in [1 - \delta, 1]} |\eta(t)| < \infty) = 1.
\end{equation}

\vspace{3mm}

 Further, suppose in addition that
\begin{equation} \label{second cond}
\exists     \, \delta \in (0,1) \ : \ \ \lim_{t \to 1^-} {\bf E} \sup_{t \in [1 - \delta, 1]} [ \ |\eta(t)|/(1 +|\eta(t))| \ ] = 0.
\end{equation}

\vspace{3mm}

 \begin{proposition}
Both the conditions \eqref{first cond} and \eqref{second cond} common are necessary and sufficient for the convergence a.e. in \eqref{what conditions}.
 \end{proposition}

\vspace{3mm}

\begin{remark}
{\rm
Evidently, if the random process $ \ \eta = \eta(t), \ t \in [0,1] $,  is such that $ \ \eta(1) = 0 $,
is continuous \emph{w.p.o.}, then the condition \eqref{what conditions} is satisfied.
}
\end{remark}

\vspace{3mm}

\hspace{3mm} Note that the exact exponential tail estimations of the maximum distribution for the {\it continuous} random fields are obtained, e.g., in \cite{Ermakov etc. 1986,Ostrovsky1}, as well as for the {\it discontinuous} random fields are in detail investigated in \cite{Ostrovsky8}.

\vspace{4mm}

Our aim in this report is to extend previous results in this direction and
   to give quantitative estimates of the norm and tail for the random variables involved in the right-hand side of \eqref{regulator}, based on ones in the left-hand side, and to extrapolate these estimates into the Grand Lebesgue Spaces of random variables.

\vspace{3mm}

\begin{remark}
{\rm Note that the systematic investigation of the convergence \emph{a.e.} in many different scientific directions, as operator and ergodic theory, functional analysis, physics, etc, is treated in \cite{Bellow Proceedings Illinois}. \par
 }
 \end{remark}

 \vspace{5mm}

\section{Lebesgue - Riesz norm estimations.}

\vspace{1mm}

Let $ (\Omega = \{\omega\}, \,  {\bf B}, \, {\bf P} ) $ be a probability space.
Recall that the classical Lebesgue-Riesz $L_p$-norm ($p\geq 1)$ of a numerical-valued random variable $\xi = \xi(\omega)$  is defined by
\begin{equation*}
||\xi||_p=||\xi||_{p,\Omega}=({\bf
E}|\xi|^p)^{1/p}=\left(\int_\Omega|\xi(\omega)|^p {\bf
P}(d\omega)\right)^{1/p}, \ \ \ p\geq 1,
\end{equation*}
\begin{equation*}
  ||\xi||_{\infty}  := \esssup_{\omega \in \Omega} |\xi(\omega)|, \ \ p=\infty.
\end{equation*}

\vspace{2mm}

Let us mention a previous result (see, e.g., \cite[Theorem 3.1]{Ostrovsky5}). Let $ \ \{ \kappa_n \}, \ n \in \mathbb N $, be a sequence of random variables such that, for some {\it single} value $ \ p \in (1,\infty) $,
$$
\lim_{n \to \infty} ||\kappa_n||_p = 0.
$$
 \ Then a factorization $ \  |\kappa_n| \le \delta_n \cdot \tau $ holds, where $ \ \delta_n \ $ is a certain positive {\it deterministic} sequence tending to zero
 and the random variable  $ \ \tau \ $  belongs to the space $ \ L_p(\Omega) $, such that
$$
\lim_{n \to \infty} \delta_n = 0, \ \ ||\tau||_p = 1.
$$

\vspace{2mm}

 \hspace{3mm} We follow further one of the auxiliary results in \cite[Lemma 2.1]{Kloeden Neuekirch}, (see also \cite[Lemma 33]{Es Sebaiy1}).

 \begin{lemma} \label{lemma Kloeden} 
 Let $\alpha>0$ and $K(p)\in[0,\infty)$ a non-negative function for $p\geq 1$. Let $Z_n$, $n\in \mathbb N$, be a sequence of random variables such that
 \begin{equation} \label{alpha cond}
\left( {\bf E}|Z_n|^p \right)^{1/p} \le K(p) \cdot n^{-\alpha},
\end{equation}
 for all $p\geq 1$ and all $n\in \mathbb N$. Then, for all $\epsilon\in (0,\alpha)$, there exists a non-negative random variable $\eta_{\alpha,\epsilon}$ such that
 \begin{equation} \label{Kloeden}
|Z_n| \le \eta_{\alpha,\epsilon} \cdot n^{-\alpha + \epsilon} \ \ \ \hbox{almost surely}
\end{equation}
for all $n\in \mathbb N$. 
 \end{lemma}

\vspace{1mm}

 \hspace{3mm} In other words,
\begin{equation} \label{overline eta}
\forall \epsilon \in (0,\alpha) \ \Rightarrow \ \eta_{\alpha,\epsilon} \stackrel{def}{=}
\sup_{n\in \mathbb N} \  \left[ \ n^{ \alpha-\epsilon} \ |Z_n|  \ \right] < \infty, \ \ \ \hbox{\emph{with probability one}}.
 \end{equation}

 \vspace{1mm}

 \ Moreover, one can conclude in addition that
\begin{equation} \label{eta epsilon}
 \exists \, p_0 = p_0(\alpha, \epsilon) \ge 1 \ : \  \forall p \ge p_0 \ \Rightarrow  \ {\bf E}|\eta_{\alpha,\epsilon}|^p < \infty.
\end{equation}

\vspace{1mm}

 Of course, only the case $ \ 0 < \epsilon < \alpha \ $ is not trivial for us.\par

 \vspace{1mm}

 Interesting applications of these results in the theory of statistic of random processes and fields may be found in \cite{Douissi}. Moreover in \cite{Kloeden Neuekirch} applications of Lemma \ref{lemma Kloeden} are given for It\^{o} stochastic differential equations.

\vspace{1mm}

  \ We want to clarify slightly the last result, under the same conditions, and extend one under additional assumptions.

\vspace{1mm}

\begin{theorem}\label{theorem eta}
 Let $\alpha>0$ and $K(p)\in[0,\infty)$ a non-negative function for $p\geq 1$. Let $Z_n$, $n\in \mathbb N$, be a sequence of r.v. satisfying the condition \eqref{alpha cond}. In addition let $ p$ and $ \epsilon $  be such that
 \begin{equation} \label{restric}
\epsilon \in (0,\min\{1,\alpha\}), \  \ \ p > 1/\epsilon.
 \end{equation}
 Then the r.v.
$$
 \eta_{\alpha,\epsilon}  \stackrel{def}{=} \sup_{n \ge 2} \left\{ \ n^{\alpha - \epsilon} \ |Z_n| \ \right\}
$$
satisfies
\begin{equation} \label{Leb Riesz est}
|| \eta_{\alpha,\epsilon} ||_p \le K(p) \ ( p \epsilon - 1)^{-1/p}.
\end{equation}
\end{theorem}

\vspace{1mm}

\begin{proof}
 We follow the same arguments as in \cite{Kloeden Neuekirch}. Indeed,
 using Tchebychev-Markov's inequality, for all $ \delta > 0 \ $,
 $ \ \epsilon \in (0,\alpha) \ $ and $ \ p > 1/\epsilon$, we have
$$
{\bf P} \left( \ n^{\alpha - \epsilon} \ |Z_n| > \delta \ \right) \le \frac{{\bf E} |Z_n|^p}{\delta^p} \cdot n^{(\alpha - \epsilon) p} \le
\frac{K^p(p)}{\delta^p}  \cdot n^{- p \epsilon},
$$
therefore
$$
\sum_{n=2}^{\infty} {\bf P} \left( \ n^{\alpha - \epsilon} \ |Z_n| > \delta  \ \right) < \infty.
$$
 \ It follows from Borel-Cantelli lemma
$$
 \lim_{n \to \infty} n^{\alpha - \epsilon} \ |Z_n| = 0
$$
with probability one. Consequently, the r.v. $\eta_{\alpha,\epsilon}$ is finite \emph{a.e}.  Moreover, for an arbitrary number $ \ q \ $ such that $ \ q > 1/\epsilon \ $ and using \eqref{alpha cond} and the integral test for the series, we get
\begin{equation*}
\begin{split}
||\eta_{\alpha,\epsilon}||_q^q & = {\bf E}  \left[ \eta_{\alpha,\epsilon} \right]^q  = {\bf E} \left\{ \ \sup_{n \ge 2}  \ \left[ \ n^{(\alpha - \epsilon)q} \ |Z_n|^q \ \right] \ \right\}\\
& \leq {\bf E} \left\{ \ \sum_{n \ge 2} n^{(\alpha - \epsilon)q} \ |Z_n|^q \ \right\} = \sum_{n \ge 2} {\bf E} \left\{ \ n^{(\alpha - \epsilon)q} \ |Z_n|^q \ \right\}\\
& \leq \sum_{n=2}^{\infty} n^{(\alpha - \epsilon)q} \ {\bf E} |Z_n|^q \ \le K^q(q) \sum_{n=2}^{\infty} n^{-q \epsilon}\\
&  \leq K^q(q)\int_1^\infty x^{-q\epsilon} \, dx = K^q(q) / (q \epsilon - 1).
\end{split}
\end{equation*}
  Substituting $ \ q = p, \ $ we obtain the absert.
 \end{proof}

\begin{remark}
{\rm
It is sufficient for the conclusion of Theorem \ref{theorem eta} to suppose only the finiteness of the function $ \ K = K(p) \ $
at least for a {\it single} value $ \ p \ $ such that $ \ p > 1/\epsilon. \ $
}
\end{remark}

 \vspace{5mm}

\section{Grand Lebesgue Spaces norm estimations.}

\vspace{1mm}

 \ We recall here for the reader convenience some known definitions and  facts  from the theory of Grand Lebesgue Spaces (GLS).
    \ Let $ b\in(1,\infty]$ and $ \ \psi = \psi(p), \ p \in [1,b)$, be a positive measurable numerical valued function, not necessarily finite in $b$, such that
    $$ \inf_{p \in [1,b)} \psi(p) > 0. $$
We denote with $ G \Psi(b)$ the set of all such functions
$\psi(p), \ p \in [1,b)$, for some  $ 1 < b \le \infty  $, and
put
$$
G\Psi := \displaystyle\bigcup_{1 < b \le \infty} G \Psi(b).
$$
For instance, if $m = {\const} > 0$,
$$
\psi_m(p) := p^{1/m},  \ \ p \in [1,\infty)
$$
or, for $ 1  < b < \infty, \ \alpha,\beta = \const \ge 0$,
$$
   \psi_{b; \alpha,\beta}(p) := (p-1)^{-\alpha} \ (b-p)^{-\beta}, \  \ p \in
   (1,b),
$$
the above functions belong to $ G \Psi(b)$.

The (Banach) Grand Lebesgue Space$ \ G \psi  = G\psi(b)$ consists
of all the real (or complex) numerical valued measurable functions
$f: \Omega \to \mathbb R$ having finite norm
\begin{equation} \label{norm psi}
    ||f||_{G\psi} \stackrel{def}{=} \sup_{p \in [1,b)} \left[ \frac{||f||_p}{\psi(p)} \right].
 \end{equation}
We write $ \ G\psi \ $ when $ \ b= \infty. \ $ The
function $ \  \psi = \psi(p) \  $ is named the {\it generating
function } for the space $G \psi$ . \par
If for instance
$$
  \psi(p) = \psi_{r}(p) = 1, \ \  p = r;  \ \  \ \psi_{r}(p) = +\infty,   \ \ p \ne r,
$$
 where $ \ r = {\const} \in [1,\infty),  \ C/\infty := 0, \ C \in \mathbb R, \ $ (an extremal case), then the correspondent
 $ \  G\psi^{(r)}(p)  \  $ space coincides  with the classical Lebesgue - Riesz space $ \ L_r = L_r(\Omega). \ $

 \begin{remark}\label{classical grand Lebesgue}
 {\rm
 Let $1<q<\infty$ and $b=q$. Define $\psi(p)=(q-p)^{-1/p}$, \ $p\in (1,q)$, and take $\Omega\subset \mathbb R^n$, \ $n\geq 1$, a measurable set with finite Lebesgue measure; then replacing in \eqref{norm psi} $p$ with $q-\varepsilon$, \ $\varepsilon\in(0,q-1)$, the space $G\psi$ reduces to 
 the classical Grand Lebesgue space $L^{q)}(\Omega)$ defined by the norm
%
%
 \begin{equation*}
 ||f||_{L^{q)}(\Omega)}=||f||_{q)}=\sup_{0<\varepsilon<q-1} \varepsilon^{\frac{1}{q-\varepsilon}}||f||_{q-\varepsilon}.
 \end{equation*}
 }
 \end{remark}

\vspace{4mm}

These spaces and their particular cases $L^{q)}(\Omega)$ (Remark \ref{classical grand Lebesgue}) are investigated in many works (e.g.
\cite{Fiorenza2000,KozOs,Liflyand_Ostrovsky_2010,Ostrovsky2,Ostrovsky4,Ostrovsky5,Samko-Umarkhadzhiev,
Samko-Umarkhadzhiev-addendum}). For example
they play an important role in the theory of Partial Differential Equations (PDEs) (see, e.g., \cite{Ahmed-Fiorenza-Formica-Gogatishvili-Rakotoson,Fiorenza-Formica-Gogatishvili-DEA2018,fioformicarakodie2017,Greco-Iwaniec-Sbordone-1997}), in interpolation
theory (see, e.g.,
\cite{AFF2022_JFAA,fioforgogakoparakoNA,fiokarazanalanwen2004}), in Functional Analysis (see, e.g., \cite{Fiorenza-Formica,FOS2021_Math_Nachr,FOS2021_JPDOA}), in the theory of
Probability (\cite{Ermakov etc.
1986,FOS2022_Contemporary Mathematics,ForKozOstr_Lithuanian,Formica_Ostrovsky_Hungarica,Ostrovsky3,Ostrovsky5}), in Statistics
(\cite[chapter 5]{Ostrovsky1}), in theory of random fields
(\cite{KozOs,Ostrovsky4}).

These spaces are rearrangement invariant (r.i.) Banach function
spaces; the fundamental function has been studied in
\cite{Ostrovsky4}. They not coincide, in the general case, with the
classical Banach rearrangement functional spaces: Orlicz, Lorentz,
Marcinkiewicz, etc., (see \cite{Liflyand_Ostrovsky_2010,Ostrovsky2}).

\hspace{3mm}  In the general case these spaces are non - separable. \par

\vspace{0.5mm}

The belonging of a function $ f: \Omega \to \mathbb{R}$ to some $ G\psi$
space is closely related to its tail function behavior
$$
 T_f(t) \stackrel{def}{=} {\bf P}(|f| \ge t), \ \ t \ge 0,
 $$
as $ \ t \to 0+ \ $ as well as when $ \ t \to \infty $ (see
\cite{KozOs,KozOsSir2019,ForKozOstr_Lithuanian}).


 \hspace{3mm} In detail, let $  \upsilon  $ be a non-zero r.v. belonging to some Grand Lebesgue Space $ \ G\psi \ $ and suppose  $ \ ||\upsilon||_{G\psi} = 1 $. Define the Young-Fenchel (or Legendre) transform of the function $\ h(p)=h[\psi](p)=p \ln \psi(p) \ $, 
 (see, e.g., \cite{Buld Koz AMS} and references therein)
$$
h^*(v)= h^*[\psi](v) := \sup_{p \in \Dom[\psi]}(p v - p \ln \psi(p))
$$
where $ \ \Dom[\psi] \ $ denotes the domain of definition (and finiteness) for the function $ \ \psi(\cdot). \ $

\vspace{1mm}
\ We get
\begin{equation} \label{Young Fen}
T_{\upsilon}(t) \le \exp(-h^*(\ln t)), \ \ t \ge e.
\end{equation}
 \ Note that the inverse conclusion, under suitable natural conditions, is true (\cite{KozOs,KozOsSir2019}). \par
 Notice also that these spaces coincide, up to equivalence of the norms and under appropriate conditions, with the so-called {\it exponential Orlicz} spaces, see e.g.
\cite{KozOs,Kozachenko at all 2018,Musielak,Ostrovsky1,Ostrovsky3}.

 \vspace{2mm}

 \begin{remark}
 {\rm
 Suppose that for some r.v. $ \ \zeta \ $ is satisfied the following condition
 $$
 \exists \, a\in (1,b) \ : \ \nu(r)\stackrel{def}{=} ||\zeta||_r   < \infty, \  \ \forall r \ge a.
 $$
 From Lyapunov's inequality follows
 $$
||\zeta||_p \le \nu(a), \ \ \forall p \in [1, a] .
$$
Therefore this r.v. $ \ \zeta \ $ belongs to the Grand Lebesgue space $ \ G\theta \ $ with the generating function $ \ \theta(\cdot) \ $ such that
$$
\theta(p) = \nu(a), \ \ p \in [1,a];  \ \ \ \theta(p) = \nu(p), \ \ p\in(a,b).
$$
Let now the r.v. $ \ \zeta  \ $ such that $||\zeta||_b < \infty, \ b\in(1,\infty]$.

\vspace{1mm}

The so-called {\it natural function} for $\zeta$  is the function $ \ \psi(p) = \psi[\zeta](p) \ $ defined by
$$
\psi(p)=\psi[\zeta](p) \stackrel{def}{=} ||\zeta||_p,
$$
with the corresponding domain of definition $ \ {\Dom} \ \psi[\zeta](\cdot). \ $

\vspace{2mm}

Analogously, let $ \ \cal{V} \ $ be an arbitrary set and let $ \  \{\zeta_v \}, \ v \in \cal{V}, \ $  be a {\it family } of random variables such that
\begin{equation} \label{family}
\sup_{v \in \cal{V}} ||\zeta_v||_p < \infty, \ \ \forall p \in [1,b) , \ \ b\in(1,\infty] \ .
\end{equation}
The {\it natural function} for this family is defined by
\begin{equation} \label{nat fun}
\psi[ \zeta_v](p) \stackrel{def}{=} \sup_{v\in \cal{V}} ||\zeta_v||_p,
\end{equation}
with the corresponding domain of definition $ \ {\Dom} \ \psi[\zeta_v]. \ $ \par
Evidently,
$$
\sup_{v \in \cal{V}} ||\zeta_v||_{G\psi[\zeta_v]} = 1.
$$

 }
 \end{remark}

\vspace{5mm}

\hspace{3mm} {Let us return to the initial relation} \eqref{alpha cond}, which we rewrite as follows: let $\alpha>0$ and $ \epsilon \in (0, \min\{1,\alpha\}) $; assume that there exists a generating function  $ \psi(p)$, $ p\in (1/\epsilon, b)$, $ 1/\epsilon < b \le \infty $  such that

\begin{equation} \label{alpha  super cond}
\left( {\bf E}|Z_n|^p \right)^{1/p} \le \psi(p) \cdot n^{-\alpha}, \hspace{4mm}  1/\epsilon < p < b.
\end{equation}

\vspace{3mm}

 \ Note that the last relation may be rewritten on the language of $ \ G\psi \ $ spaces as follows
\begin{equation} \label{alpha  GPSI cond}
||Z_n||_{G\psi} \le  n^{-\alpha}.
\end{equation}

 \vspace{3mm}

 \ Define as before the  {\it key} r.v.
\begin{equation*}
 \eta_{\epsilon} \stackrel{def}{=} \sup_{n \ge 1} \left\{ \ n^{\alpha - \epsilon} \ |Z_n| \ \right\},
\end{equation*}
which is finite almost everywhere.  Moreover, by virtue of Theorem \ref{theorem eta}, we have
\begin{equation} \label{psi est}
 || \ \eta_{\epsilon} \ ||_p  \le  \psi(p) \ ( p \epsilon - 1)^{-1/p}, \ \ \ p \in (1/\epsilon,b) \ .
\end{equation}
 Introduce the following new generating function
\begin{equation}\label{new generating}
\kappa_{\alpha, \epsilon}[\psi](p) \stackrel{def}{=} \psi(p) \ ( p \epsilon - 1)^{-1/p},  \ \ \ p \in (1/\epsilon,b),
\end{equation}
which define the corresponding Grand Lebesgue space $ \ G \kappa_{\alpha, \epsilon}[\psi]$.

\vspace{3mm}

 \ To summarize, we have the following result. \par

 \vspace{4mm}

\begin{theorem}\label{theorem estimate in GLS}
 Let $\alpha>0$, \ $b\in(1,\infty]$ and  $ \  \epsilon \in (0, \min\{1,\alpha\}) $. Let $Z_n$, $n\in \mathbb N$, be a sequence of r.v. and let $ \ \psi(p)$,  $ \ p\in (1/\epsilon, b)$, be a function satisfying the condition \eqref{alpha  super cond}. Let $\kappa_{\alpha, \epsilon}[\psi]$ the function defined in \eqref{new generating}.
 Then the r.v. $\eta_{\epsilon}$ defined above
satisfies
\begin{equation} \label{Gpsi ineq}
||  \eta_{\epsilon}  ||_{G \kappa_{\alpha, \epsilon}[\psi]} \le  1,
\end{equation}
and the corresponding {\it exponential} tail estimate \eqref{Young Fen}. \par

\end{theorem}
\vspace{3mm}

\begin{remark}
{\rm
 Note that one can also consider the space $ \  G \psi(a, b), \  1 < a < b \le \infty, \ $
 where $ \  p \in (a, b). \ $ In this case the norm is defined alike above

\vspace{3mm}

\begin{equation} \label{norm psi  ab}
  ||f||_{G\psi(a,b)}  \stackrel{def}{=} \sup_{p \in (a,b)} \left[ \frac{||f||_p}{\psi(p)} \right].
\end{equation}
}
\end{remark}
\vspace{5mm}

\section{A slight generalization.}

\vspace{1mm}

 \hspace{3mm} Consider the source random sequence, which
 we denote now as  $ \ \{V_n\}, \ n \in \mathbb N $,
 taking values in a certain (complete) Banach space $ \ X  \ $ equipped with the norm $ \ \langle \cdot \rangle = \langle \cdot \rangle_X, \ $
 non necessarily separable, and let us impose the following condition

\begin{equation} \label{new restriction}
|| \ \langle V_n \rangle_X||_p \le \psi(p) \cdot \epsilon_n,
\end{equation}
where  $ \  \psi = \psi(p)$, $ \ p \in [1,b),  \ b \in (1,\infty]$, is some positive  generating function and
 $ \  \{ \epsilon_n \} \ $ is a numerical positive sequence tending to zero as
 $ \ n \to \infty: \ \displaystyle \lim_{n \to \infty} \epsilon_n = 0. \ $ \par
 \ Moreover, we suppose that

\begin{equation} \label{p condit}
\exists \ p \in [1,b)  \ : \ \sum_{n=1}^{\infty} \epsilon_n^p  < \infty.
\end{equation}

\vspace{3mm}

 \ Let also $ \ \{\beta_n\} \ $ be another numerical positive sequence  also tending to zero and such that
 the  following function

 \begin{equation} \label{sigma function}
 \sigma^p(p) \stackrel{def}{=} \sum_{n=1}^{\infty} \frac{\epsilon^p_n}{\beta^p_n},
 \end{equation}
 is finite at least for some value $ \ p \in [1,b). \ $ It is easy to ground the existence of this function under the condition \eqref{p condit}. \par

\vspace{3mm}

\begin{example}
{\rm
Let here

$$
 \epsilon_n = n^{-\alpha} \ \ln^m(n+ 1), \ \ \alpha > 0, \ \ m \in \mathbb R;
 $$
then one can  choose

$$
\beta_n := n^{-\theta} \ \ln^{-\nu}(n+1),  \ \ \theta \in (0,\alpha), \ \ \nu \in \mathbb R,
$$
and correspondingly one can choose  $ p $ an arbitrary number such that $ \ p > 1/(\alpha - \theta) $,  for arbitrary values of the parameters $ \ m, \nu; \ $
but in the case when

$$
\frac{m + \nu}{\alpha - \theta} < - 1
$$
one can choose  even  $ \  p := 1/(\alpha - \theta). \ $ \par

\vspace{3mm}

 Of course, one can use arbitrary {\it slowly varying} at infinity functions, strictly  positive, $  L =L(n), \ M = M(n) $, instead of $ \ \ln^k(\cdot). \ $ \par
 \ In detail, let
$$
 \epsilon_n := n^{-\alpha} \ L(n), \ \ \alpha > 0,
 $$
then one can  analogously  choose

$$
\beta_n := n^{-\theta} \ M(n),  \ \ \theta \in (0,\alpha),
$$
with $ \ p > 1/(\alpha - \theta). \ $ \par

}

\end{example}

\vspace{4mm}

 \ Put also
\begin{equation} \label{eta rv}
\zeta \stackrel{def}{=} \sup_n \left[ \ \frac{\langle V_n\rangle_X}{\beta_n} \ \right].
\end{equation}

\vspace{3mm}

 \ We deduce quite alike in the proof of Theorem \ref{theorem eta}

\vspace{3mm}

$$
{\bf E} \zeta^p \le \psi^p(p) \sum_{n = 1}^{\infty} \left[ \ \frac{\epsilon_n^p}{\beta_n^p} \ \right] = \psi^p(p) \cdot  \sigma^p(p).
$$

\vspace{3mm}

 \hspace{3mm} To summarize: \par

\begin{theorem}\label{theorem generalization}
Under the above restrictions and notations the following estimate holds

\begin{equation}  \label{Main general}
||\zeta||_{p, \Omega}  \le \psi(p) \cdot  \sigma(p),
\end{equation}
or equally, if we introduce the generating function $ \ \gamma(p) :=  \psi(p) \cdot  \sigma(p),  \ $ in the interval of its finiteness, then the norm of $\ \zeta$ in the Grand Lebesgue space $ \ G\gamma \ $ satisfies
\begin{equation}  \label{Main simpl}
||\zeta||_{G{\gamma}}  \le 1.
\end{equation}

\end{theorem}

\vspace{3mm}

\begin{example}
{\rm

 Let now $ \ \epsilon_n = q^n, \ n = 0,1,2,\ldots; \ q ={\const} \in (0,1). \ $ One can choose
 $ Q = {\const} \in (q,1), \ $  so that $ \  \delta := q/Q \in (0,1). \ $   We get

$$
\sigma^p(p) = \sum_{n=0}^{\infty} \frac{q^{p n}}{Q^{p  n}} = \sum_{n=0}^{\infty} \delta^{ \ p n \ } = \frac{1}{1 - \delta^p} < \infty
$$
for all the {\it positive} values $ \ p. \ $  But {\it for the values}  $ \  p \ge 1  \ $ we  have

$$
\sigma(p) = (1 - \delta^p)^{ -1/p} \le (1- \delta)^{-1},
$$
and we conclude, on the basis of Theorem \ref{theorem generalization},

\begin{equation}  \label{Main non improvable}
||\zeta||_{p, \Omega}  \le \psi(p)/(1 - \delta).
\end{equation}
}
\end{example}

\vspace{5mm}

\section{Examples. Lower estimates.}

\vspace{1mm}

 \hspace{3mm}   \ {\bf 1.} \ Let us ground that the estimate of Theorem \ref{theorem eta} is  essentially non-improvable,
 of course, under our restrictions. Namely, let  $ \ \alpha = {\const} > 0 $ and put

$$
Z_n := \frac{\theta_n}{n^{\alpha}}, \  \ \ n = 2,3,\ldots,
$$
where $ \  \{\theta_n \} \ $ is a sequence of {\it independent} positive identical distributed r.v. having  standard exponential
distribution
$$
{\bf P}(\theta_n > t) = e^{-t}, \  \ \ t > 0.
$$

\vspace{3mm}

\hspace{3mm} Note first of all that, since the r.v. $ \ \{Z_n\} \ $ are independent, the convergence $ \ Z_n \to 0$ as  $ n \to \infty $
is completely equivalent to  the following restriction
$$
\forall \epsilon > 0 \ \Rightarrow \ \sum_{n = 1}^{\infty} {\bf P} (Z_n> \epsilon)  < \infty,
$$
or equally
$$
\forall \epsilon > 0 \ \Rightarrow \ \sum_{n = 1}^{\infty}  \exp(- \epsilon \ n^{\alpha})   < \infty.
$$
 \ Evidently,  both these conditions are satisfied, therefore really

$$
{\bf P}  ( \ \lim_{n \to \infty} Z_n \to 0) = 1,
$$
and consequently $ \ \sup_n Z_n < \infty \ $ a.e.

\vspace{3mm}

 \ Define, as before, $ \ \displaystyle \eta := \sup_{n \ge 1} \left[ \ n^{\alpha - \epsilon} \ Z_n \ \right], \ \epsilon \in (0,\alpha). \ $ It is not hard
to prove that as $ \ u \to \infty, \ (u \ge 3)$,
\begin{equation} \label{examp exp estim}
{\bf P}(\eta > u) \asymp C(\epsilon)  \ u^{-1/\epsilon}, \ \ \ \ C(\epsilon)= \frac{1}{\epsilon} \Gamma\left(\frac{1}{\epsilon}\right) = \Gamma\left(1+\frac{1}{\epsilon}\right),
\end{equation}
where $\Gamma(\cdot)$ is the Euler's Gamma function.
\vspace{1mm}

 \ In detail, denoting $ \ \kappa_n = \theta_n/n^{\epsilon}=n^{\alpha-\epsilon}Z_n $, by integral test for the series and the change of variable $z=u x^{\epsilon}$, as $ \ u \to \infty$ we have
 \begin{equation*}
 \begin{split}
 \Sigma_1 \stackrel{def}{=} \ &\sum_n {\bf P}(\kappa_n > u) = \sum_n \ \exp \left(- u \ n^{\epsilon} \right) \ \sim \int_0^{\infty} \exp \left( \ - u \ x^{\epsilon}  \ \right) \ dx\\
 & = u^{-{1}/{\epsilon}} \cdot \frac{1}{\epsilon} \int_0^\infty z^{\frac{1}{\epsilon}-1}e^{-z}\, dz \ \asymp \ \frac{1}{\epsilon} \Gamma\left(\frac{1}{\epsilon}\right) u^{-{1}/{\epsilon}} \ =  C(\epsilon) \ u^{-1/\epsilon}.
 \end{split}
 \end{equation*}

From the Bonferroni's inequality (upper bound), as $ \ u \to \infty$, we have
\begin{equation*}
{\bf P}(\eta > u)  \ = {\bf P}(\sup_n \kappa_n > u) \ ={\bf P} \left(\bigcup_n \kappa_n > u \right) \le \Sigma_1.
\end{equation*}
On the other hand, taking into account the independence of the random variables $ \ \{ \kappa_n\}  $, we have
\begin{equation*}
\begin{split}
 \Sigma_2  & \stackrel{def}{=}  \ \sum_n \sum_{m > n} \ {\bf P}(\ \kappa_n > u, \ \kappa_m > u) = \sum_n \sum_{m > n} \ {\bf P}(\ \kappa_n > u) \ {\bf P}(\ \kappa_m > u) \\
& \leq \sum_n \sum_m  \ {\bf P}(\ \kappa_n > u) \ {\bf P}(\ \kappa_m > u)  = \Sigma_1^2.
\end{split}
\end{equation*}

 \hspace{3mm} The relation \eqref{examp exp estim} follows from the \emph{bilateral} Bonferroni's estimate and $\Sigma_2\leq \Sigma_1^2$ since, for $u\to\infty$, we have
$$
\Sigma_1 \asymp \Sigma_1 - \Sigma_1^2 \ \leq \ \Sigma_1 - \Sigma_2 \ \le \ {\bf P}(\sup_n \kappa_n > u) \ \le \ \Sigma_1.
$$

\vspace{2mm}

 \ Further, we have therefore

\begin{equation} \label{Lower bound first}
||\eta||_p \asymp (1/\epsilon - p)^{-1/p}, \ p \in [1, \ 1/\epsilon).
\end{equation}

\vspace{3mm}

 \ \hspace{3mm} \ {\bf 2.} \ On the other hand, $ \ \eta \ge Z_1, \ $ therefore also

\vspace{3mm}

\begin{equation} \label{Lower bound second}
 ||\eta||_p \ge ||Z_1||_p \ge C \, p, \ \ \ p \ge 1.
\end{equation}

 Thus, in this case

\begin{equation} \label{Lower bound third}
 ||\eta||_p \ge  C \max \{ \ p, \ (1/\epsilon - p)^{-1/p} \ \}, \ \ \ p \in [1, \ 1/\epsilon).
\end{equation}

\vspace{3mm}

 \ \hspace{3mm} \ {\bf 3.} \ Let the inequality \eqref{psi est} be given under conditions of the first section.
 Let us assume  also that the  generating function $ \ \psi = \psi(p) \ $  is the natural function for the variable $ \ Z_1: \ $
$$
\psi(p) = ||Z_1||_p, \ \ \ p \in [1,b).
$$

\vspace{2mm}

 \hspace{3mm} Suppose now that the $ \, p \, $ belongs to the interval $ \ (a,b), \ $ where $ \ a \ge 4/\epsilon, \ b \in (a,\infty]. \ $
 \ Then one has for such values $ \ p \ $


\begin{equation} \label{D estimate}
||\eta||_p \le 3^{1/\epsilon} \cdot \psi(p).
\end{equation}


 On the other hand,
\begin{equation} \label{unit case}
||\eta||_p \ge ||Z_1||_p = \psi(p).
\end{equation}

\vspace{4mm}

 \section{Concluding remarks.}

\vspace{1mm}

 \hspace{4mm} It is interesting in our opinion to investigate a more general convergence of random variables almost everywhere in the metric spaces, even in the topological spaces, in the spirit of  Remark \ref{another criterion}.

\vspace{6mm}

\emph{\textbf{\footnotesize Acknowledgements}.} {\footnotesize M.R. Formica is member of Gruppo
Nazionale per l'Analisi Matematica, la Probabilit\`{a} e le loro Applicazioni (GNAMPA) of the Istituto Nazionale di Alta Matematica (INdAM) and member of the UMI group \lq\lq Teoria dell'Approssimazione e Applicazioni (T.A.A.)\rq\rq and is partially supported by the INdAM-GNAMPA project, {\it Risultati di regolarit\`{a} per PDEs in spazi di funzione non-standard}, codice CUP\_E53C22001930001 and partially supported by University of Naples \lq\lq Parthenope\rq\rq, Dept. of Economic and Legal Studies, project CoRNDiS, DM MUR 737/2021, CUP I55F21003620001.
}

\vspace{2mm}

\textbf{Author contributions:} All authors contribute equally to the manuscript.

\vspace{2mm}

\textbf{Conflict of interest}: The authors declare that they have no conflict of interest.

\vspace{2mm}

\textbf{Data Availability}: Data sharing is not applicable to this article as no datasets were generated or analyzed during the current study.

%

\end{document}